\documentclass[11pt]{amsart}
\usepackage{amsmath,amsfonts,amsthm,amscd,amssymb,graphicx}

\newtheorem{theorem}{Theorem}[section]
\newtheorem{lemma}[theorem]{Lemma}
\newtheorem{proposition}[theorem]{Proposition}
\newtheorem{corollary}[theorem]{Corollary}

\newtheorem{remark}[theorem]{Remark}

\newcommand{\ip}{\int_\mathbb{R}}

\def\bU{{\hat{u}}}
\def\bV{{\hat{v}}}
\def\R{\Re e}
\def\I{\Im m}

\begin{document}
\title[Stability of viscous shocks in isentropic gas dynamics]{Stability of viscous shocks in isentropic\\ gas dynamics}
\author[Barker, Humpherys, Rudd, and Zumbrun]{Blake Barker, Jeffrey Humpherys,\\ Keith Rudd, and Kevin Zumbrun}

\date{Last Updated:  March 23, 2007}

\thanks{ This work was supported in part by the National Science Foundation award numbers DMS-0607721 and DMS-0300487.}

\address{Department of Mathematics, Brigham Young University, Provo, UT 84602}
\email{jeffh@math.byu.edu}
\address{Department of Mathematics, Indiana University, Bloomington, IN 47402}
\email{kzumbrun@indiana.edu}

\maketitle

\begin{abstract}
In this paper, we examine the stability 
problem for 
viscous shock solutions of
the isentropic compressible Navier--Stokes equations,
or $p$-system with real viscosity.   
We first revisit the work of Matsumura and Nishihara, extending the known parameter regime for which small-amplitude viscous shocks are provably spectrally stable by an optimized version of their original argument. Next, using a novel spectral energy estimate, we show that there are no purely real unstable eigenvalues for any shock strength.

By related estimates, we show that unstable eigenvalues are confined to a bounded region independent of shock strength.  Then through an extensive numerical Evans function study, we show that there is no unstable spectrum in the entire right-half plane, thus demonstrating numerically that large-amplitude shocks are spectrally stable up to Mach number $M\approx 3000$ for $1 \le \gamma \leq 3$.  This strongly suggests that shocks are stable independent of amplitude and the adiabatic constant $\gamma$.  We complete our study by showing that finite-difference simulations of perturbed large-amplitude shocks converge to a translate of the original shock wave, as expected.
\end{abstract}

\section{Introduction}
\label{int}

Consider the isentropic compressible Navier-Stokes equations in one 
spatial dimension
expressed in Lagrangian coordinates,
also known as the $p$-system with real viscosity:
\begin{equation}
\begin{split}
\label{psystem}
v_{t}-u_{x} &=0, \\
u_{t}+p(v)_{x} &= \left(\frac{u_x}{v}\right)_x,
\end{split}
\end{equation}
where, physically, $v$ is the specific volume, 
$u$ is the velocity,
and $p(v)$ is the pressure law.   We assume that $p(v)$ is adiabatic, that is, a $\gamma$-law gas satisfying $p(v) = a_0 v^{-\gamma}$ for some constants $a_0>0$ and $\gamma \geq 1$.  In the thermodynamical rarified gas approximation, $\gamma>1$ is the average over constituent particles of $\gamma=(n+2)/n$, where $n$ is the number of internal degrees of freedom of an individual particle \cite{Ba}: $n=3$ ($\gamma=1.66...$) for monatomic, $n=5$ ($\gamma=1.4$) for diatomic gas.  More generally, $\gamma$ is usually taken within $1 \leq \gamma \leq 3$ in physical approximations of gas- or fluid-dynamical flow \cite{Sm}.

This system is an important and widely studied gas-dynamical model (see for example \cite{Sm} and references within), and yet little is presently known about the stability of its large-amplitude viscous shock wave solutions.  Over two decades ago, Matsumura and Nishihara \cite{MN} used a clever energy estimate to show that small-amplitude shocks are stable under zero-mass perturbations.  The linear portion of their analysis is equivalent to proving spectral stability, which through the more recent work of Zumbrun and collaborators \cite{ZH,MZ.1,MZ.2,MZ.3,Z.2,HZ}, implies asymptotic orbital stability, hereafter called nonlinear stability.  We remark that the results of \cite{MZ.2,MZ.3,Z.2} hold
for shocks of arbitrary amplitude, and thus nonlinear stability is implied by spectral stability.  Hence, for large-amplitude shocks, spectral stability is the missing piece of the stability puzzle and the subject of our present focus.

In this paper, we first improve upon the work in \cite{MN} slightly by extending the known parameter regime for which small-amplitude viscous shocks are provably spectrally stable, using an optimized version of the same method.  We also show that this method cannot be extended any further to larger amplitudes.  Using a novel spectral energy estimate, however, we are able to show that there are no purely real unstable eigenvalues for any shock strength.  We note that this result is stronger than that which could be given by the Evans function stability index (sometimes called the orientation index), which only measures the parity of unstable real eigenvalues, see \cite{GJ,BSZ,GZ,Z.2}.  A consequence of this result (which follows also by the Evans function computations used to determine the stability index \cite{Z.2}) is that if an instability were to occur for large-amplitude viscous shocks, its onset, or indeed any change in the number of unstable eigenvalues, would be associated with a Hopf-like bifurcation in which one or more conjugate pairs of eigenvalues cross through the imaginary axis; see \cite{TZ.1,TZ.2,TZ.3} for further discussion of this scenario.  

Continuing our investigations, we appeal to numerical Evans function computation to explore the large-amplitude regime through the use of winding number calculations via the argument principle.  Before doing so, however, we rule out high-frequency instability through the combination of two spectral energy estimates, showing that unstable eigenvalues are confined to a bounded region independent of shock amplitude.  This reduces the problem to investigation, feasible by numerics, of a compact set.  Then by checking the low-frequency regime by repeating several Evans function computations, we determine whether or not a particular viscous shock is spectrally stable.   As a final verification, we use a finite-difference method to simulate perturbed large-amplitude viscous shocks, and check whether they converge, as expected, to a translate of the original profile.

{\bf Conclusions and results of numerical investigations:}
We carry out our numerical experiments far into the hypersonic shock regime,  exploring up through Mach number $M\approx 3000$ for $1 \le \gamma \leq 3$.  Particular attention is given to the monatomic and diatomic cases, $\gamma=5/3$ and $\gamma=7/4$, respectively.  In all cases, our results are {\it consistent with spectral stability}  (hence also linear and nonlinear stability \cite{MZ.2,MZ.3,Z.2}).  This strongly suggests that viscous shock profiles in an isentropic are spectrally stable independent of both amplitude and $\gamma$.  Our bounds on the size of unstable eigenvalues, which are independent of shock strength, may be viewed as a first step in establishing such a result analytically.

{\bf Extensions and open questions.}
The present study is not a numerical proof.  However, as discussed in \cite{Br.1}, it could be converted to 
numerical proof by the implementation of interval arithmetic and a posteriori error estimates for numerical solution of ODE. This would be an interesting direction for future investigation.  A crucial step in carrying out numerical proof by interval arithmetic is by analytical estimates special to the problem at hand to sufficiently reduce the computational domain to make the computations feasible in realistic time.
This we have carried out by the surprisingly strong estimates of Section \ref{hfsec} and Appendices \ref{profbds}--\ref{initbds}.

A second very interesting direction would be to establish stability in the large-amplitude limit by a singular-perturbation analysis, an avenue we intend to follow in future work. Together, these two projects would give a complete, rigorous proof of stability for arbitrary-amplitude viscous shock solutions of \eqref{psystem} on the physical range $\gamma\in [1,3]$.

\section{Background}

By a viscous shock profile of \eqref{psystem}, we mean a traveling wave solution
\begin{equation*}
\begin{split}
v(x,t)&=\bV(x-s t),\\
u(x,t)&=\bU(x-s t),
\end{split}
\end{equation*}
moving with speed $s$ and having asymptotically constant end-states $(v_{\pm},u_{\pm})$. As an alternative, we can translate $x \rightarrow x-s t$, and consider instead stationary solutions of
\begin{equation}
\begin{split}
\label{movingpsystem}
v_t - s v_x - u_x &= 0,\\
u_t - s u_x + (a_0 v^{-\gamma})_x &= \left(\frac{u_x}{v}\right)_x.
\end{split}
\end{equation}
Under the rescaling 
$(x,t,v,u) \rightarrow 
(-\varepsilon sx, \varepsilon s^2 t, v/\varepsilon, -u/(\varepsilon s))$, 
where $\varepsilon$ is chosen so that $0 < v_+ < v_- = 1$,  
our system takes the form
\begin{equation}
\begin{split}
\label{rescaled}
v_t + v_x - u_x &= 0,\\
u_t + u_x + (a v^{-\gamma})_x &= \left(\frac{u_x}{v}\right)_x,
\end{split}
\end{equation}
where $a = a_0 \varepsilon^{-\gamma-1} s^{-2}$.  Thus, the shock profiles of \eqref{rescaled} satisfy the ordinary differential equation
\begin{equation*}
\begin{split}
v' - u' &= 0,\\
u' + (a v^{-\gamma})' &= \left(\frac{u'}{v}\right)',
\end{split}
\end{equation*}
subject to the boundary conditions $(v(\pm \infty),u(\pm \infty)) =(v_{\pm},u_{\pm})$. This simplifies to
\begin{equation*}
v' + (a v^{-\gamma})' = \left(\frac{v'}{v}\right)'.
\end{equation*}
By integrating from $-\infty$ to $x$, we get the profile equation
\begin{equation}
\label{profile}
v' = v(v-1 + a  (v^{-\gamma}-1)),
\end{equation}
where $a$ is found by setting $x=+\infty$, thus yielding the Rankine-Hugoniot condition
\begin{equation}
\label{RH}
a = -\frac{v_+ - 1}{v_+^{-\gamma} - 1} = v_+^\gamma \frac{1-v_+}{1-v_+^\gamma}.
\end{equation}
Evidently, $a\to \gamma^{-1}$
in the weak shock limit $v_+\to 1$, while 
$ a\sim v_+^\gamma $ in the strong shock limit $v_+\to 0$.

\begin{remark}
Since the profile equation \eqref{profile} is first order scalar, 
it has a monotone solution.  
Since $v_+<v_-$, we have that ${\bV}_x<0$ for all $x\in \mathbb{R}$.
\end{remark}

By linearizing \eqref{rescaled} about the profile $(\bV,\bU)$, we get the eigenvalue problem
\begin{equation}
\label{eigen1}
\begin{split}
&\lambda v + v' - u' =0,\\
&\lambda u + u' - \left(\frac{h(\bV)}{\bV^{\gamma+1}}v\right)' = \left(\frac{u'}{\bV}\right)',
\end{split}
\end{equation}
where
\begin{equation}
\label{f}
h(\bV) = -\bV^{\gamma+1} + a(\gamma-1) + (a+1) \bV^\gamma
\end{equation}
We say that a shock profile of \eqref{psystem} is spectrally stable if the linearized system \eqref{eigen1}  has no spectrum in the closed deleted right half-plane
\[
P = \{\R(\lambda) \geq 0\}\setminus\{0\},
\]
that is, there are no growth or oscillatory modes for \eqref{eigen1}.  We remark that a traveling wave profile always has a zero-eigenvalue associated with its translational invariance.  This generally negates the possibility of good uniform bounds in energy estimates, and so we employ the standard technique (see \cite{Go.1,ZH}) of transforming into integrated coordinates. This goes as follows:

Suppose that $(v,u)$ is an eigenfunction of \eqref{eigen1} with a non-zero eigenvalue $\lambda$.  Then
\[
\tilde{u}(x) = \int_{-\infty}^x u(z) dz , \quad
\tilde{v}(x) = \int_{-\infty}^x v(z) dz,
\]
and their derivatives decay exponentially as $x \rightarrow \infty$. Thus, by substituting and then integrating, $(\tilde{u},\tilde{v})$ satisfies (suppressing the tilde)
\begin{subequations}\label{ep}
\begin{align}
&\lambda v + v' - u' =0, \label{ep:1}\\
&\lambda u + u' -  \frac{h(\bV)}{\bV^{\gamma+1}} v' = \frac{u''}{\bV}.\label{ep:2}
\end{align}
\end{subequations}
This new eigenvalue problem differs spectrally from \eqref{eigen1} only at $\lambda=0$, hence spectral stability of \eqref{eigen1} is implied by spectral stability of \eqref{ep}.  Moreover, since \eqref{ep} has no eigenvalue at $\lambda=0$, one can expect to have a better chance of developing a successful spectral energy method to prove stability.  We demonstrate this in the following section.

\section{Stability of small-amplitude shocks}\label{small}

In this section we further the work in \cite{MN} by extending slightly the known parameter regime for which small-amplitude viscous shocks are provably spectrally stable.  We also show that this method cannot be extended any further for larger amplitudes. 

\begin{theorem}[\cite{MN}]
\label{smallamp}
Viscous shocks of \eqref{psystem} are spectrally stable whenever
\begin{equation}
\label{condition}
\left(\frac{v_+^{\gamma+1}}{a\gamma}\right)^2+2(\gamma-1)\left(\frac{v_+^{\gamma+1}}{a\gamma}\right)-(\gamma-1) \geq 0
\end{equation}
In particular, as $v_+\rightarrow 1$ (hence $a\gamma\rightarrow 1$), the left-hand side of \eqref{condition} approaches $\gamma$ and so the inequality is satisfied.  Therefore, small-amplitude viscous shocks of \eqref{psystem} are spectrally stable.
\end{theorem}

\begin{proof}

We note that $h(\bV) > 0$.  By multiplying \eqref{ep:2} by both the conjugate $\bar{u}$ and $\bV^{\gamma+1}/h(\bV)$ and integrating along $x$ from $\infty$ to $-\infty$, we have
\[
\ip \frac{\lambda u \bar{u}\bV^{\gamma+1}}{h(\bV)} + \ip \frac{u' \bar{u}\bV^{\gamma+1}}{h(\bV)} -  \ip v' \bar{u} = \ip \frac{u''\bar{u}\bV^\gamma}{h(\bV)}.
\]
Integrating the last three terms by parts and appropriately using \eqref{ep:1} to substitute for $u'$ in the third term gives us
\[
\ip \frac{\lambda |u|^2 \bV^{\gamma+1}}{h(\bV)} + \ip \frac{u' \bar{u}\bV^{\gamma+1}}{h(\bV)} + \ip v (\overline{\lambda v + v'}) + \ip \frac{\bV^\gamma|u'|^2}{h(\bV)} = -\ip \left(\frac{\bV^\gamma}{h(\bV)}\right)' u'\bar{u}.
\]
We take the real part and appropriately integrate by parts to get
\[
\R(\lambda)\ip \left[ \frac{\bV^{\gamma+1}}{h(\bV)}|u|^2+|v|^2 \right] +  \ip g(\bV) |u|^2 + \ip \frac{\bV^\gamma}{h(\bV)}|u'|^2= 0,
\]
where
\[
g(\bV) = - \frac{1}{2} \left[\left(\frac{\bV^{\gamma+1}}{h(\bV)}\right)' + \left(\frac{\bV^\gamma}{h(\bV)}\right)'' \right].
\]
Thus, to prove stability it suffices to show that $g(\bV) \geq 0$ on $[v_+,1]$.

By straightforward computation, we obtain identities:
\begin{align}
\gamma h(\bV) - \bV h'(\bV) &= a\gamma(\gamma-1) + \bV^{\gamma+1}\quad\mbox{and}\label{l1}\\
\bV^{\gamma-1}\bV_x &= a\gamma - h(\bV)\label{I2}.
\end{align}
Using \eqref{l1} and \eqref{I2}, we abbreviate a few intermediate steps below:
\begin{align}
g(\bV) &= -\frac{\bV_x}{2}\left[ \frac{(\gamma+1)\bV^\gamma h(\bV) - \bV^{\gamma+1}h'(\bV)}{h(\bV)^2} + \frac{d}{d\bV}\left[ \frac{\gamma \bV^{\gamma-1}h(\bV)-\bV^\gamma h'(\bV)}{h(\bV)^2} \bV_x \right]\right]\notag\\
&= -\frac{\bV_x}{2}\left[ \frac{\bV^\gamma\left((\gamma+1)h(\bV) - \bV h'(\bV)\right)}{h(\bV)^2} + \frac{d}{d\bV}\left[ \frac{\gamma h(\bV)-\bV h'(\bV)}{h(\bV)^2} (a\gamma-h(\bV)) \right]\right]\notag\\
%
&=-\frac{a\bV_x\bV^{\gamma-1}}{2 h(\bV)^3} \times\notag\\
& \qquad\left[ \gamma^2(\gamma+1)\bV^{\gamma+2} - 2 (a+1)\gamma(\gamma^2-1)\bV^{\gamma+1}+(a+1)^2\gamma^2(\gamma-1)\bV^\gamma\right.\notag\\
&\qquad\qquad +\left. a\gamma(\gamma+2)(\gamma^2-1)\bV-a (a+1) \gamma^2 (\gamma^2-1) \right]\notag\\
&=-\frac{a\bV_x\bV^{\gamma-1}}{2 h(\bV)^3}[(\gamma+1)\bV^{\gamma+2}+\bV^\gamma(\gamma-1)\left((\gamma+1)\bV-(a+1)\gamma\right)^2 \label{preMN}\\
&\qquad + a\gamma(\gamma^2-1)(\gamma+2)\bV-a (a+1)\gamma^2(\gamma^2-1)]\notag\\
&\geq -\frac{a\bV_x\bV^{\gamma-1}}{2 h(\bV)^3}[(\gamma+1)\bV^{\gamma+2}+ a\gamma(\gamma^2-1)(\gamma+2)\bV-a (a+1)\gamma^2(\gamma^2-1)]\notag\\
&\geq-\frac{\gamma^2 a^3 \bV_x (\gamma+1)}{2 h(\bV)^3 v_+}\left[\left(\frac{v_+^{\gamma+1}}{a\gamma}\right)^2+2(\gamma-1)\left(\frac{v_+^{\gamma+1}}{a\gamma}\right)-(\gamma-1)\right].\label{MN}
\end{align}
Thus from \eqref{condition}, we have spectral stability.\qed
\end{proof}

\begin{figure}[t]
\begin{center}
$\begin{array}{cc}
\includegraphics[width=5.75cm]{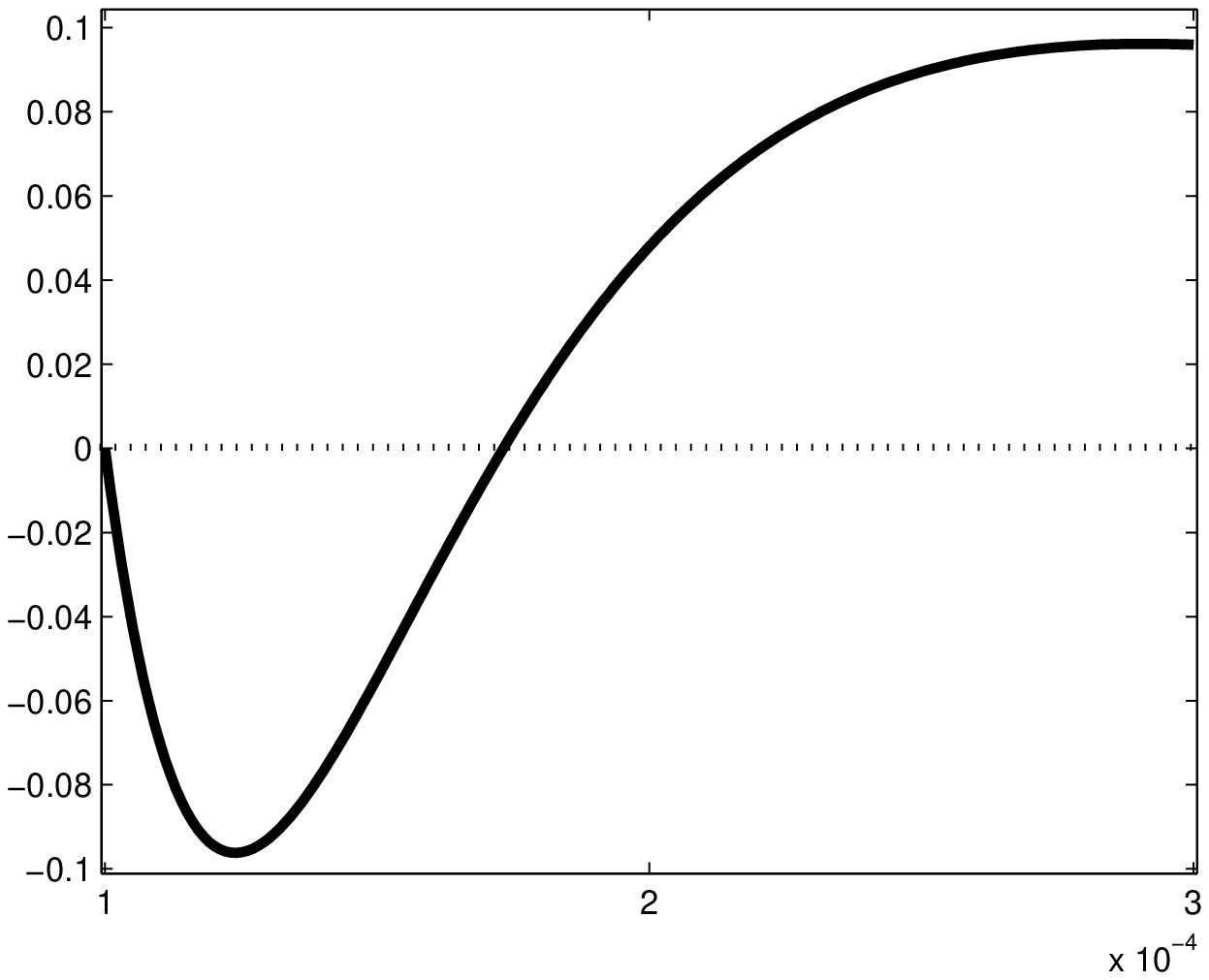} & \includegraphics[width=5.75cm]{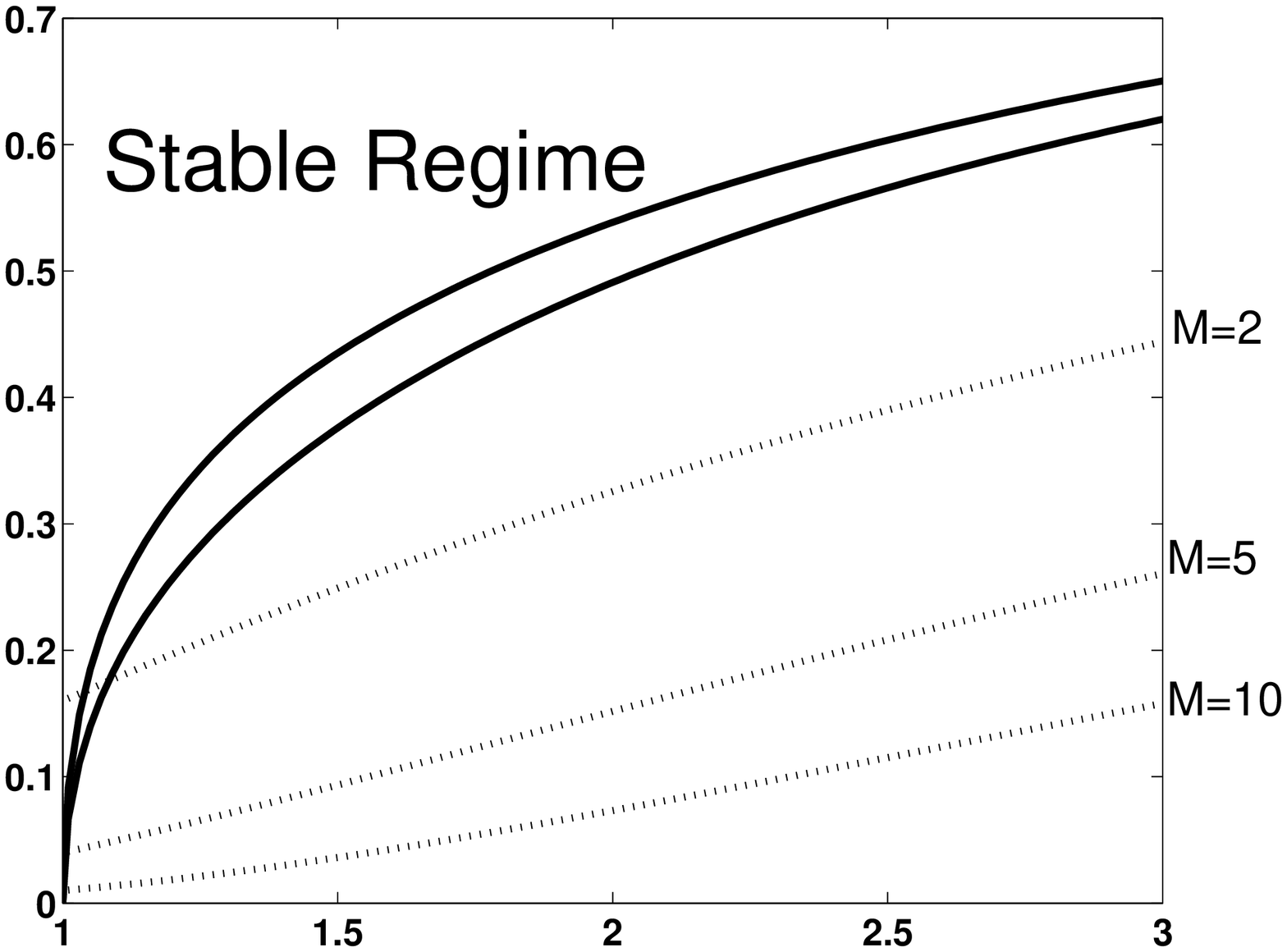} \\
\mbox{\bf (a)} & \mbox{\bf (b)}
\end{array}$
\end{center}
\caption{
In (a), we have a graph of $g(\hat v)$ against $\hat v$ for $v_+= 1\times 10^{-4}$ and $\gamma=2.0$.  Note that $g(\bV)$ dips down below zero on the left-hand size.  Hence, the energy estimate will not generalize beyond the small-amplitude regime.  In (b) we graph the stability boundaries (dark lines) given by \eqref{condition} and \eqref{condition3}, where $\gamma$ and $v_+$ are the horizontal and vertical axes, respectively.  We see that in our scaling \eqref{condition3} is only a modest improvement over \eqref{condition}.  The dotted lines correspond from top to bottom as the parameter regimes for Mach numbers $2$, $5$, and $10$ (see the appendix to see how to determine the Mach number).  Hence, this energy estimate does not even hold for shocks at Mach 2 and $\gamma > 1.084$.}
\label{dip_g}
\end{figure}

We note that  the hypothesis in \eqref{condition} is not sharp.  Indeed, one can show from \eqref{preMN}, that a stronger condition could be given as
\begin{align}
&(\gamma+1)\bV^{\gamma+2}+\bV^\gamma(\gamma-1)\left((\gamma+1)\bV-(a+1)\gamma\right)^2 \label{condition2} \\
&\quad + a\gamma(\gamma^2-1)(\gamma+2)\bV-a (a+1)\gamma^2(\gamma^2-1) \geq 0,\notag
\end{align}
which is sharp in the following sense:  When this condition fails to be true, then $g(\bV)$ is no longer nonnegative, and thus the energy method fails.  In Figure \ref{dip_g}(a), we see that $g(\bV)$ dips on the left-hand side when this inequality is compromised.  We remark further that near $v_+$, the left-hand side of \eqref{condition2} is monotone increasing in $\hat v$.
Thus,
\eqref{condition2} holds if and only if
\begin{align}
&(\gamma+1)v_+^{\gamma+2}+v_+^\gamma(\gamma-1)\left((\gamma+1)v_+-(a+1)\gamma\right)^2 \label{condition3} \\
&\quad + a\gamma(\gamma^2-1)(\gamma+2)v_+-a (a+1)\gamma^2(\gamma^2-1) \geq 0.\notag
\end{align}
We see from Figure \ref{dip_g}(b) that \eqref{condition3} is only a marginal improvement over \eqref{condition}.  However, since \eqref{condition3} is sharp, we cannot hope to prove large-amplitude spectral stability using this approach.  Instead, we proceed by a combined analytical and numerical approach as in \cite{Br.1,Br.2,BrZ}, first showing that unstable eigenvalues can occur only in a bounded
set, then searching for eigenvalues in this set by computing the Evans function numerically.  Before doing so, however, we show in the following section that real unstable eigenvalues do not exist, even for large-amplitude viscous shocks.

\section{No real unstable eigenvalues}

In this section we use a novel spectral energy estimate to show that there are no purely real unstable eigenvalues for any shock strength.  We note that this result is stronger than that which could be given by the Evans function stability index (sometimes called the orientation index), which only measures the parity of unstable real eigenvalues, see \cite{BSZ,GZ,Z.2}.  The fact that this holds for all shock strengths is interesting because it is among the strongest statements about large-amplitude spectral stability that has been proven to date.

\begin{theorem}
Viscous shocks of \eqref{psystem} have no unstable real spectra.
\end{theorem}
\begin{proof}
We multiply \eqref{ep:2} by the conjugate $\bar{v}$ and integrate along $x$ from $\infty$ to $-\infty$.  This gives
\[
\ip \lambda u \bar{v} + \ip u' \bar{v} - \ip \frac{h(\bV)v' \bar{v}}{\bV^{\gamma+1}} = \ip \frac{u''\bar{v}}{\bV}.
\]
Notice that on the real line, $\bar{\lambda}=\lambda$.  Using \eqref{ep:1} to substitute for $\overline{\lambda v}$ in the first term and for $u''$ in the last term, we get
\[
\ip u (\bar{u}'-\bar{v}') + \ip u' \bar{v} - \ip \frac{h(\bV)v' \bar{v}}{\bV^{\gamma+1}} = \ip \frac{(\lambda v' + v'')\bar{v}}{\bV}.
\]
Separating terms and simplifying gives
\[
\ip u \bar{u}' + 2 \ip u' \bar{v} - \ip \frac{h(\bV)v' \bar{v}}{\bV^{\gamma+1}}  = \lambda  \ip \frac{v'\bar{v}}{\bV} + \ip \frac{v'' \bar{v}}{\bV}.
\]
We further simplify by substituting for $u'$ in the second term and integrating  the last terms by parts to give,
\[
\ip u \bar{u}' + 2 \ip (\lambda v + v') \bar{v} - \ip \frac{h(\bV)v' \bar{v} }{\bV^{\gamma+1}} = \lambda  \ip \frac{v'\bar{v}}{\bV} + \ip \frac{\bV_x}{\bV^2}v' \bar{v} - \ip  \frac{|v'|^2}{\bV},
\]
which yields
\[
\ip u \bar{u}' + 2 \lambda \ip |v|^2 + 2 \ip v' v -  a\gamma \ip  \frac{v' \bar{v}}{\bV^{\gamma+1}}  + \ip \frac{|v|^2}{\bV} = \lambda  \ip \frac{v'\bar{v}}{\bV}.
\]
By taking the real part, we arrive at
\[
\frac{\lambda}{2} \ip \left(4 - \frac{\bV_x}{\bV^2} \right) |v|^2 - \frac{a\gamma(\gamma+1)}{2} \ip  \frac{\bV_x}{\bV^{\gamma+2}} |v|^2 + \ip \frac{|v'|^2}{\bV} = 0.
\]
This is a contradiction when $\lambda \geq 0$.\qed
\end{proof}
We remark that the absence of positive real eigenvalues limits the admissible onset of instability to Hopf-like bifurcations where a pair of conjugate eigenvalues crosses the imaginary axis.  In the following section, we give an upper bound on the spectral frequencies that are admissible.  In other words, we show that if a pair of conjugate eigenvalues cross the imaginary axis, they must do so within these bounds.

\section{High-frequency bounds}\label{hfsec}

In this section, we prove high-frequency spectral bounds.  This is an important step because it gives a ceiling as to how far along both the imaginary and real axes that one must explore for spectrum when doing Evans function computations.  Indeed if we want to check for roots of the Evans function in the unstable half-plane, say using the argument principle, then we need only compute within these bounds.  If no roots are found therein, then we have  strong numerical evidence that the given shock is spectrally stable.  (Indeed, at the expense of further effort, such a calculation may be used as the basis of numerical proof, as described in \cite{Br.1,Br.2}.)  In this section we show that the high-frequency bounds are quite strong, only allowing unstable eigenvalues to persist in a relatively small triangle adjoining the origin.  Moreover, these bounds are {\it independent of the shock amplitude}.

\begin{lemma}
The following identity holds for $\R \lambda \geq 0$:
\begin{align}
(\R(\lambda) + |\I (\lambda)|) & \ip \bV |u|^2 - \frac{1}{2}\ip \bV_x |u|^2 + \ip |u'|^2\notag\\
 &\leq \sqrt{2} \ip \frac{h(\bV)}{\bV^\gamma} |v'| |u| +  \ip \bV |u'||u|\label{id1}.
\end{align}
\end{lemma}

\begin{proof}
We multiply \eqref{ep:2} by $\bV {\bar u}$ and integrate along $x$.  This yields
\[
\lambda \ip \bV |u|^2 + \ip \bV u'\bar{u} + \ip |u'|^2 = \ip \frac{h(\bV)}{\bV^\gamma} v'\bar{u} .
\]
We get \eqref{id1} by taking the real and imaginary parts and adding them together, and noting that $|\R(z)| + |\I(z)| \leq \sqrt{2}|z|$.\qed
\end{proof}

\begin{lemma}
\label{kawashima}
The following identity holds for $\R \lambda \geq 0$:
\begin{equation}
\label{id3}
\ip |u'|^2 = 2\R(\lambda)^2\ip|v|^2 + \R(\lambda)\ip \frac{|v'|^2}{\bV} + \frac{1}{2} \ip \left[\frac{h(\bV)}{\bV^{\gamma+1}} + \frac{a\gamma}{\bV^{\gamma+1}} \right] |v'|^2
\end{equation}
\end{lemma}

\begin{proof}
We multiply \eqref{ep:2} by ${\bar v'}$ and integrate along $x$.  This yields
\[
\lambda \ip u\bar{v}' + \ip u'\bar{v}' - \ip \frac{h(\bV)}{\bV^{\gamma+1}}|v'|^2 = \ip \frac{1}{\bV}u''\bar{v}' = \ip \frac{1}{\bV}(\lambda v' + v''){\bar v'}.
\]
Using \eqref{ep:1} on the right-hand side, integrating by parts, and taking the real part gives
\[
\R \left[ \lambda \ip u\bar{v}' + \ip u'\bar{v}'\right] = \ip \left[\frac{h(\bV)}{\bV^{\gamma+1}} + \frac{\bV_x}{2 \bV^2} \right] |v'|^2 + \R(\lambda)\ip \frac{|v'|^2}{\bV}.
\]
The right hand side can be rewritten as
\begin{equation}
\label{id3_1}
\R \left[ \lambda \ip u\bar{v}' + \ip u'\bar{v}'\right] = \frac{1}{2} \ip \left[\frac{h(\bV)}{\bV^{\gamma+1}} + \frac{a\gamma}{\bV^{\gamma+1}} \right] |v'|^2 + \R(\lambda)\ip \frac{|v'|^2}{\bV}.
\end{equation}
Now we manipulate the left-hand side.  Note that
\begin{align*}
\lambda \ip u\bar{v}' + \ip u'\bar{v}' &= (\lambda+\bar{\lambda}) \ip u\bar{v}' - \ip u(\bar{\lambda}\bar{v}' + \bar{v}'')\\
&= -2\R(\lambda) \ip u' \bar{v} - \ip u \bar{u}''\\
&= -2\R(\lambda) \ip (\lambda v + v') \bar{v} + \ip |u'|^2.
\end{align*}
Hence, by taking the real part we get
\[
\R \left[ \lambda \ip u\bar{v}' + \ip u'\bar{v}'\right] = \ip |u'|^2 - 2\R(\lambda)^2 \ip |v|^2.
\]
This combines with \eqref{id3_1} to give \eqref{id3}.\qed
\end{proof}

\begin{remark}
Lemma \ref{kawashima} is a special case of the high-frequency bounds given in \cite{MZ.3,Z.2}.  We also note that \eqref{id3} follows from a ``Kawashima-type'' estimate as described in \cite{MZ.3,Z.2}.
\end{remark}

\begin{lemma}\label{Hlem}
For $h(\bV)$ as in \eqref{f}, we have
\begin{equation}
\label{id2}
\sup_{\bV} \left| \frac{h(\bV)}{\bV^\gamma}\right| = \gamma 
\frac{1-v_+}{1-v_+^\gamma}
\leq \gamma.
\end{equation}
\end{lemma}

\begin{proof}
Defining
\begin{equation}\label{H}
H(\bV):=h(\bV)\bV^{-\gamma} = -\bV + a(\gamma-1)\bV^{-\gamma} + (a+1),
\end{equation}
we have $H'(\bV)= -1 -a\gamma(\gamma-1)\bV^{-\gamma-1}<0$ for $0<v_+\le \bV\le v_-= 1$,
hence the maximum of $H$ on $\bV\in [v_+,v_-]$ is achieved at
$\bV=v_+$.
Substituting \eqref{RH} into \eqref{H} and simplifying
yields \eqref{id2}.
\end{proof}

We complete this section by proving our high-frequency bounds.

\begin{theorem}
Any eigenvalue $\lambda$ of \eqref{ep} with nonnegative real part satisfies
\begin{equation}
\label{hfbounds}
\R(\lambda) + |\I(\lambda)| \leq (\sqrt{\gamma}+\frac{1}{2})^2.
\end{equation}
\end{theorem}

\begin{proof}
Using Young's inequality twice on right-hand side of \eqref{id1} together with \eqref{id2}, we get
\begin{align*}
(\R(\lambda) +& |\I (\lambda)|) \ip \bV |u|^2 - \frac{1}{2}\ip \bV_x |u|^2 + \ip |u'|^2 \\
&\leq \sqrt{2} \ip \frac{h(\bV)}{\bV^\gamma} |v'| |u| +  \ip \bV |u'||u|\\
&\leq \theta \ip \frac{h(\bV)}{\bV^{\gamma+1}} |v'|^2 + \frac{(\sqrt{2})^2}{4\theta} \ip \frac{h(\bV)}{\bV^\gamma} \bV |u|^2 + \epsilon \ip \bV |u'|^2 + \frac{1}{4 \epsilon} \ip \bV |u|^2\\
&< \theta \ip \frac{h(\bV)}{\bV^{\gamma+1}} |v'|^2  + \epsilon \ip |u'|^2 + \left[\frac{\gamma}{2\theta} + \frac{1}{4 \epsilon}\right] \ip \bV |u|^2.
\end{align*}
Assuming that $0<\epsilon<1$ and $\theta = (1-\epsilon)/2$, this simplifies to
\begin{align*}
(\R(\lambda) + |\I (\lambda)|) & \ip \bV |u|^2 + (1-\epsilon) \ip |u'|^2 \\
&<\frac{1-\epsilon}{2} \ip \frac{h(\bV)}{\bV^{\gamma+1}} |v'|^2 +  \left[\frac{\gamma}{2\theta} + \frac{1}{4 \epsilon}\right] \ip \bV |u|^2.
\end{align*}
Applying \eqref{id3} yields
\[
(\R(\lambda) + |\I (\lambda)|) \ip \bV |u|^2  <  \left[\frac{\gamma}{1-\epsilon} + \frac{1}{4 \epsilon}\right]  \ip \bV |u|^2,
\]
or equivalently,
\[
(\R(\lambda) + |\I (\lambda)|) <  \frac{(4 \gamma-1)\epsilon - 1}{4\epsilon(1-\epsilon)}.
\]
Setting $\epsilon = 1/(2\sqrt{\gamma}+1)$ gives \eqref{hfbounds}.\qed
\end{proof}

In the following section, we do an extensive Evans function study to explore numerically the rest of the right-half plane in an effort to locate the presence of unstable complex eigenvalues.

\section{Evans function computation}

In this section, we numerically compute the Evans function to locate any unstable eigenvalues, if they exist, in our system.  The Evans function $D(\lambda)$ is analytic to the right of the essential spectrum and is defined as the Wronskian of decaying solutions of the eigenvalue equation for the linearized operator \eqref{ep} (see \cite{E.1,E.2,E.3,E.4,EF,AGJ}).  In a spirit similar to the characteristic polynomial, we have that $D(\lambda)=0$ if and only if $\lambda$ is in the point spectrum of the linearized operator \eqref{ep}.  While the Evans function is generally too complex to compute explicitly, it can readily be computed numerically, even for large systems \cite{HuZ.2}.

Since the Evans function is analytic in the region of interest, we can numerically compute the winding number in the right-half plane.  This allows us to systematically locate roots (and hence eigenvalues) within.  As a result, spectral stability can be determined, and in the case of instability, one can produce bifurcation diagrams to illustrate and observe its onset.  This approach was first used by Evans and Feroe \cite{EF} and has been advanced further in various directions (see for example \cite{PW,PSW,AS,Br.1,Br.2,BDG,HuZ.2}).

We begin by writing \eqref{ep} as a first-order system $W' = A(x,\lambda) W$, where
\begin{equation}
\label{evans_ode}
A(x,\lambda) = \begin{pmatrix}0 & \lambda & 1\\0 & 0 & 1\\ \lambda \bV& \lambda\bV &f(\bV)-\lambda \end{pmatrix},\quad W = \begin{pmatrix} u\\v\\v'\end{pmatrix},\quad \prime = \frac{d}{dx},
\end{equation}
and $f(\bV) = \bV- \bV^{-\gamma} h(\bV)$, with $h$ as in \eqref{f}.  Note that eigenvalues of \eqref{ep} correspond to nontrivial solutions of $W(x)$ for which the boundary conditions $W(\pm\infty)=0$ are satisfied.  We remark that since $\bV$ is asymptotically constant in $x$, then so is $A(x,\lambda)$.  Thus at each end-state, we have the constant-coefficient system
\begin{equation}
\label{apm}
W' = A^\pm(\lambda) W.
\end{equation}
Hence solutions that satisfy the needed boundary condition must emerge from the unstable manifold $W_1^-(x)$ at $x=-\infty$ and the stable manifold $W_2^+(x) \wedge W_3^+(x)$ at $x=\infty$.  In other words, eigenvalues of \eqref{ep} correspond to the values of $\lambda$ for which these two manifolds intersect, or more precisely, when $D(\lambda): = \det(W_1^- W_2^+ W_3^+)_{\mid x=0}$ is zero.

As an alternative, we consider the adjoint formulation of the Evans function \cite{PW,BSZ}, where instead of computing the $2$-dimensional stable manifold, we find the single trajectory $\widetilde{W}_1^+(x)$ 
coming from the unstable manifold of 
\begin{equation}\label{adjode}
\widetilde{W}' = -\widetilde{W} A(x,\lambda)
\end{equation}
at $x=\infty$.  Note that $\widetilde{W}_1^+(x)$ is biorthogonal to both $W_2^+(x)$ and  $W_3^+(x)$ since $(\widetilde{W}(x) W(x))'=0$ for all $x$ and the initial data of $\widetilde{W}$ is biorthogonal to that of $W$.  Hence, the original manifolds intersect when $\widetilde{W}_1^+$ and $W_1^-$ are biorthogonal.  Therefore, the adjoint Evans function takes the form  $D_+(\lambda): = (\widetilde{W}_1^+ W_1^-)_{\mid x=0}$.

To further improve the numerical efficiency and accuracy of the shooting scheme, we rescale $W$ and $\widetilde{W}$ to remove exponential growth/decay at infinity, and thus eliminate potential problems with stiffness.  Specifically, we let $W(x) = e^{\mu^- x} V(x)$, where $\mu^-$ is the growth rate of the unstable manifold at $x=-\infty$, and we solve instead $V'(x) = (A(x,\lambda)-\mu^- I)V(x)$.  We initialize $V(x)$ at $x=-\infty$ to be the eigenvector of $A_-(\lambda)$ corresponding to $\mu^-$.  Similarly, it is straightforward to rescale and initialize $\widetilde{W}(x)$ at $x=\infty$.

Numerically, we use a finite domain for $x$, replacing the end states $x=\pm\infty$ with $x=\pm L$, for sufficiently large $L$.  Some care needs to be taken, however, to make sure that we go out far enough to produce good results.  In Appendices \ref{profbds} and \ref{initbds}, we explore the decay rates of the profile $\bV$ and $A(x,\lambda)$ and combine our analysis with numerical convergence experiments to conclude that our domain $[-L,L]$ is sufficiently large.  Our experiments, described below, were primarily conducted using $L=12$, with relative error bounds ranging mostly between $10^{-3}$ and  $10^{-4}$.  Larger choices of $L$ were needed on the high end of the Mach scale, going up to $L=18$ in some cases, to get the relative errors down to $10^{-4}$.  In Table \ref{tablerel}, we provide a sample of relative errors in $D(\lambda)$ for large-amplitude shocks. 

\begin{table}
\begin{center}
\begin{tabular}{ccccc}
\hline
$L$ & $\gamma=1.2$ & $\gamma=1.4$ & $\gamma=1.666$ & $\gamma=1.8$ \\
\hline
8 & 1.23(-1) & 1.16(-1) & 1.08(-1) & 1.04(-1) \\
10 & 2.07(-2) & 1.46(-2) & 1.75(-2) & 1.78(-2) \\
12 & 2.00(-3) & 1.40(-3) & 9.85(-4) & 7.20(-4) \\
14& 6.90(-4) & 5.31(-4) & 4.73(-4) & 4.71(-4)\\
\hline
\end{tabular}
\caption{Relative errors in $D(\lambda)$ are computed by taking the maximum relative error for 60 contour points evaluated along the semicircle in Figure \ref{evans_out}(a).  Samples were taken for varying $L$ and $\gamma$, leaving $v_+$ fixed at $v_+=10^{-4}$ (Mach $M \approx 1669$).  We used $L=8,10,12,14,16$ and $\gamma=1.2,1.4,1.666,1.8$.  Relative errors were computed using the next value of $L$ as the baseline.}
\label{tablerel}
\end{center}
\end{table}

We remark also that in order to produce analytically varying Evans function output, the initial data. $V(-L)$ and $\widetilde{V}(L)$, must be chosen analytically.  The method of Kato \cite[pg. 99]{Kato}, also described in \cite{BrZ}, does this well by replacing the eigenvectors of \eqref{apm} with analytically defined spectral projectors (see also \cite{BDG,HSZ}).

Before we can compute the Evans function, we first need to compute the traveling wave profile.  We use both Matlab's {\tt ode45} routine, which is the adaptive fourth-order Runge-Kutta-Fehlberg method (RKF45), and Matlab's {\tt bvp4c} routine, which is an adaptive Lobatto quadrature scheme.  Both methods work well and produce essentially equivalent profiles.

\begin{figure}[t]
\begin{center}
$\begin{array}{cc}
\includegraphics[width=5.75cm]{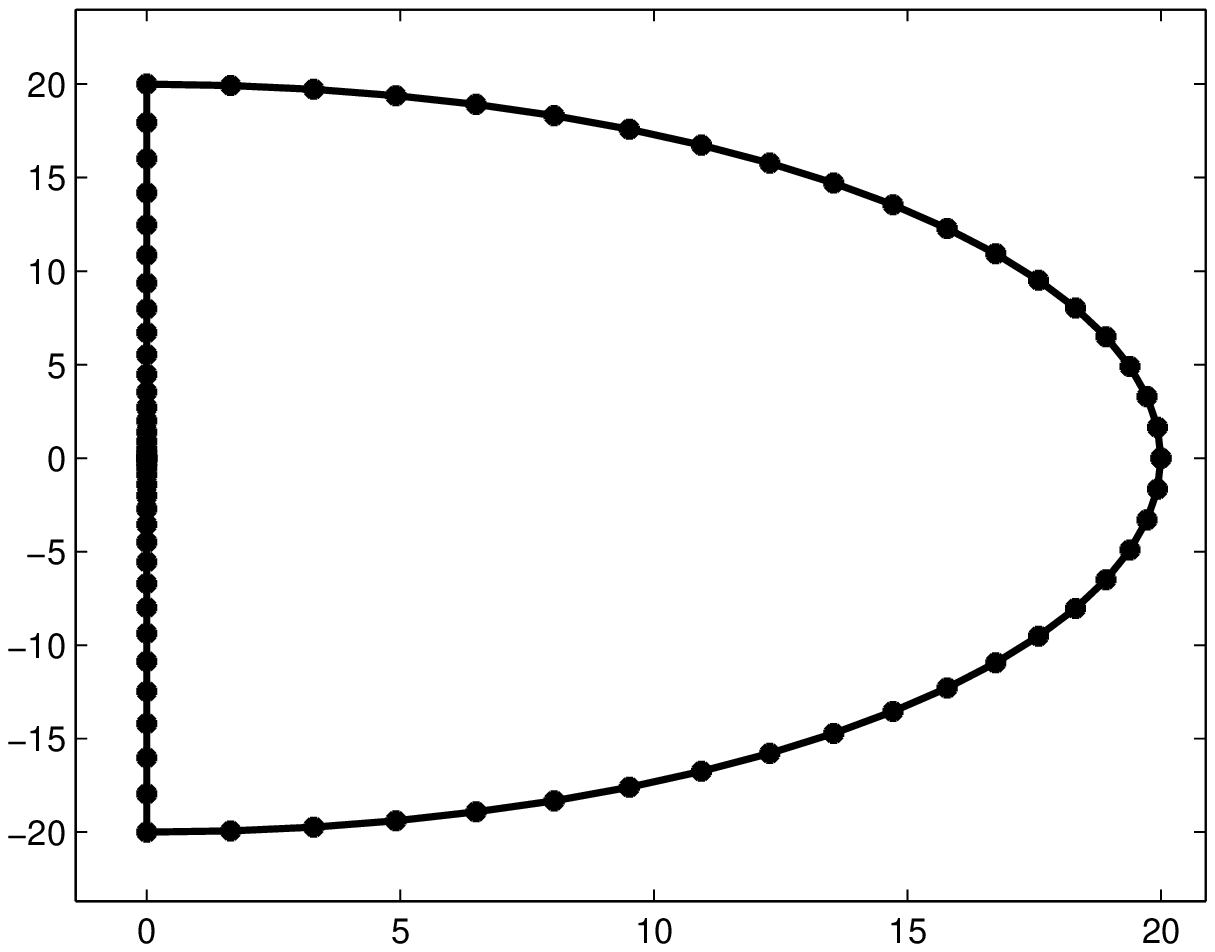} & \includegraphics[width=5.75cm]{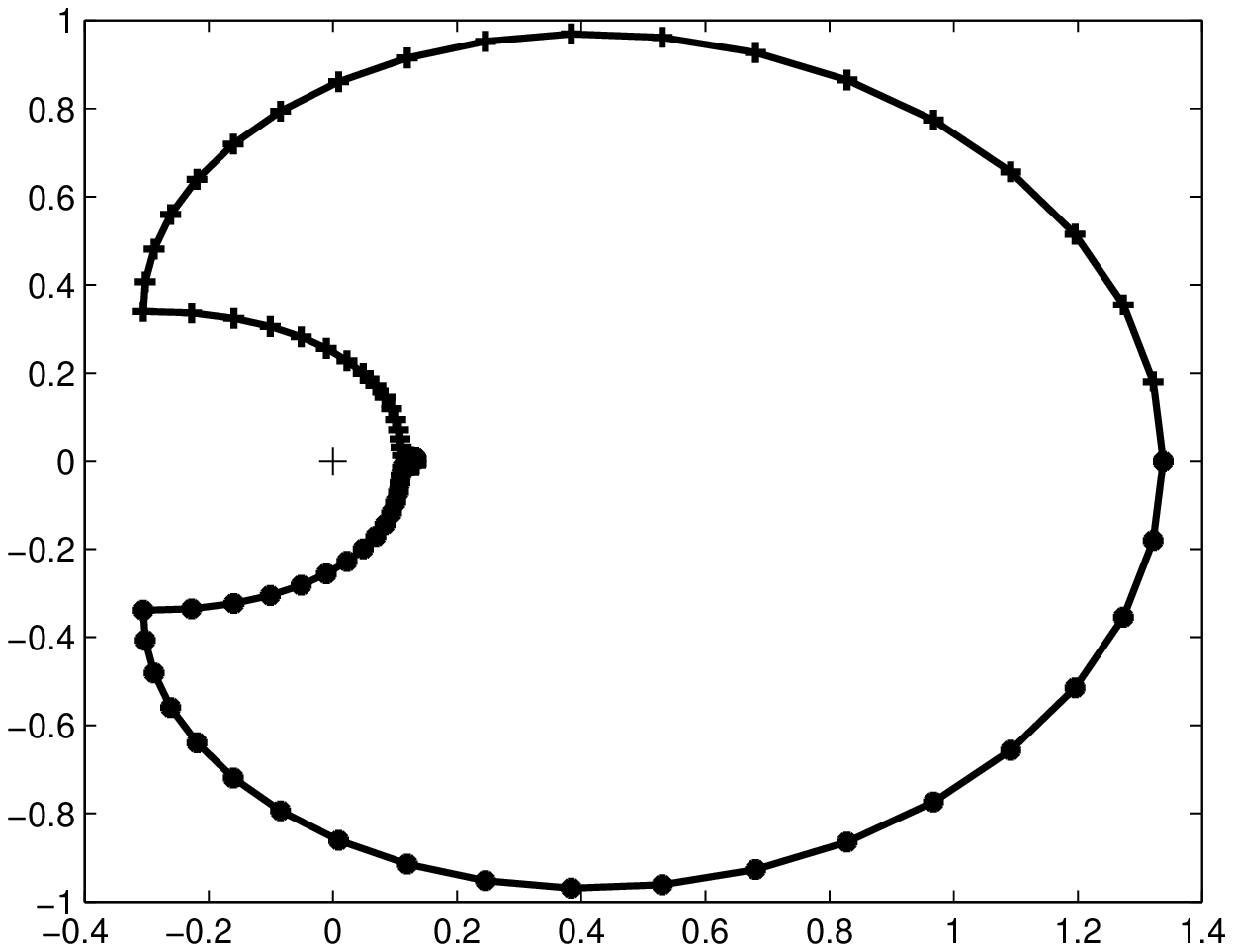} \\
\mbox{\bf (a)} & \mbox{\bf (b)}
\end{array}$
\end{center}
\caption{The graph of the contour and its mapping via the Evans function.  In (a), we have a contour, which is a large semi-circle aligned on the imaginary axis on the right-half plane (horizontal axis real, vertical axis imaginary).  
In (b) we have the image of the contour mapped by the Evans function.  Note that the winding number of this graph is zero.  Hence, there are no unstable eigenvalues in the semi-circle.  Together with the high-frequency bounds, this implies spectral stability.  Our computation was carried out for $\gamma=5/3$ and $v_+=1\times 10^{-4}$.  This corresponds to a monatomic gas with Mach number $M \approx 1669$ (see Appendix \ref{mach}).}
\label{evans_out}
\end{figure}

Our experiments were carried out on the range
\[
(\gamma, M)\in [1,3]\times [1.6,3000].
\]
Recall, for $\gamma\in [1,3]$, that shocks are known to be stable for $M \in [1,1.6]$, by \eqref{condition3}, Section \ref{small}, hence this completes the study of the range 
\[
(\gamma, M)\in [1,3]\times [1,3000]
\]
from minimum Mach number $M=1$ far into the hypersonic shock regime, and encompassing the entire physically relevant range of $\gamma$.

For each value of $\gamma$, the Mach number $M$ was varied on logarithmic scale with regular mesh from $M=1.6$ up to $M=3000$.  We used $50$ mesh points for $\gamma=1.0+0.1k$, where $k=1,2,\ldots,20$.  For the special values $\gamma=1.4$ and $1.666$ (monatomic and diatomic cases), we did a refined study with $400$ mesh points.

For each value of $(\gamma, M)$, we computed the Evans function along semi-circular contours that contain the triangular region found in the previous section via our high-frequency bounds; see Figure \ref{evans_out}(a).  The ODE calculations for individual $\lambda$ were carried out using Matlab's {\tt ode45} routine, which is the adaptive 4th-order Runge-Kutta-Fehlberg method (RKF45), and after scaling out the exponential decay rate of the constant-coefficient solution at spatial infinity, as described in \cite{Br.1,Br.2,BrZ,BDG,HuZ.2}.  This method is known to have excellent accuracy \cite{BDG,HuZ.2};  in addition, the adaptive refinement gives automatic error control.  Typical runs involved roughly $300$ mesh points, with error tolerance set to {\tt AbsTol = 1e-6} and {\tt RelTol = 1e-8}.  Values of $\lambda$ were varied on a contour with $60$ points.  As a check on winding number accuracy, it was tested a posteriori that the change in argument of $D$ for each step was less than $\pi/25$.  Recall, by Rouch\'e's Theorem, that accuracy is preserved so long as the argument varies by less than $\pi$ along each mesh interval.

In all the cases that we examined, the winding number was zero.  This indicates that the shocks we considered are spectrally stable, and in view of \cite{MZ.1,MZ.2,MZ.3}, nonlinear stability follows.  Moreover, in light of the large Mach numbers considered, this is highly suggestive of stability for all shock strengths.

\section{Numerical evolution of strong shocks}

In this section, we use a standard finite-difference method to simulate perturbed large-amplitude viscous shocks, and show that they converge, as expected, to a translate of the original profile.

We do this with a nonlinear Crank-Nicholson scheme with a Newton solver to deal with the nonlinearity.  This method provides second-order accuracy and will allow for larger time steps than a naive explicit scheme.  By differentiating the viscosity term in \eqref{rescaled}, we have

\begin{figure}[t]
\begin{center}$
\begin{array}{lr}
\includegraphics[width=5.75cm]{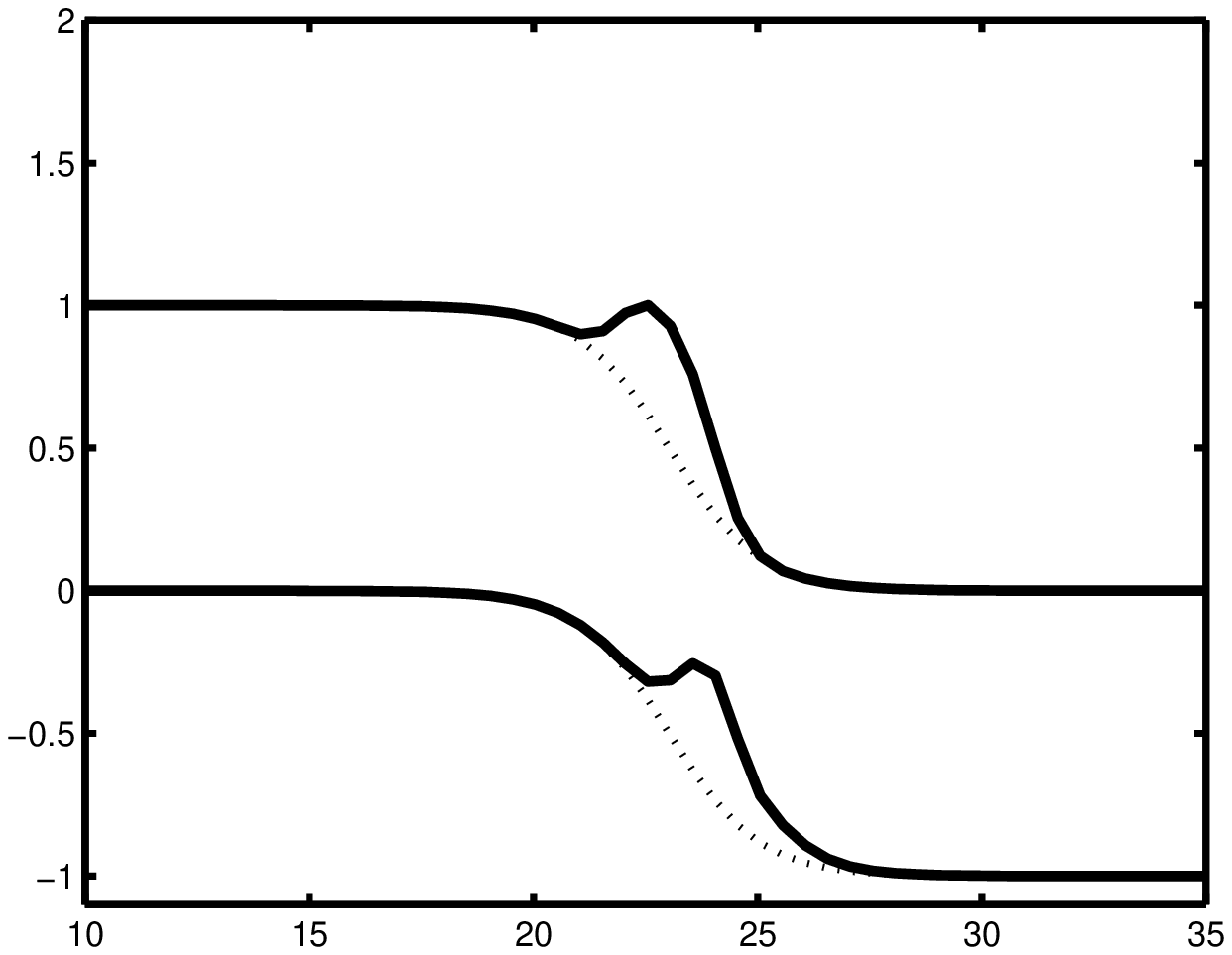} &
\includegraphics[width=5.75cm]{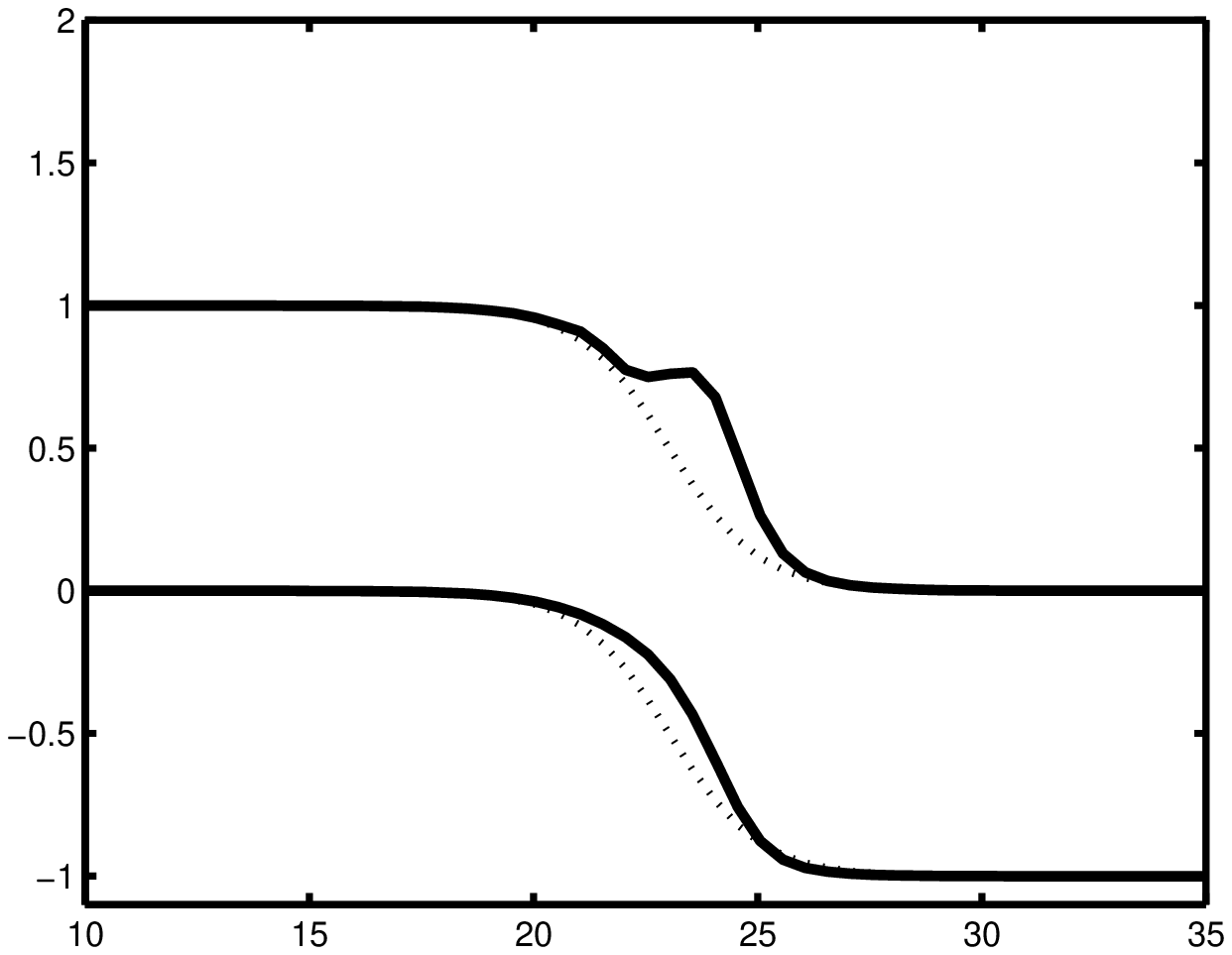} \\ 
\includegraphics[width=5.75cm]{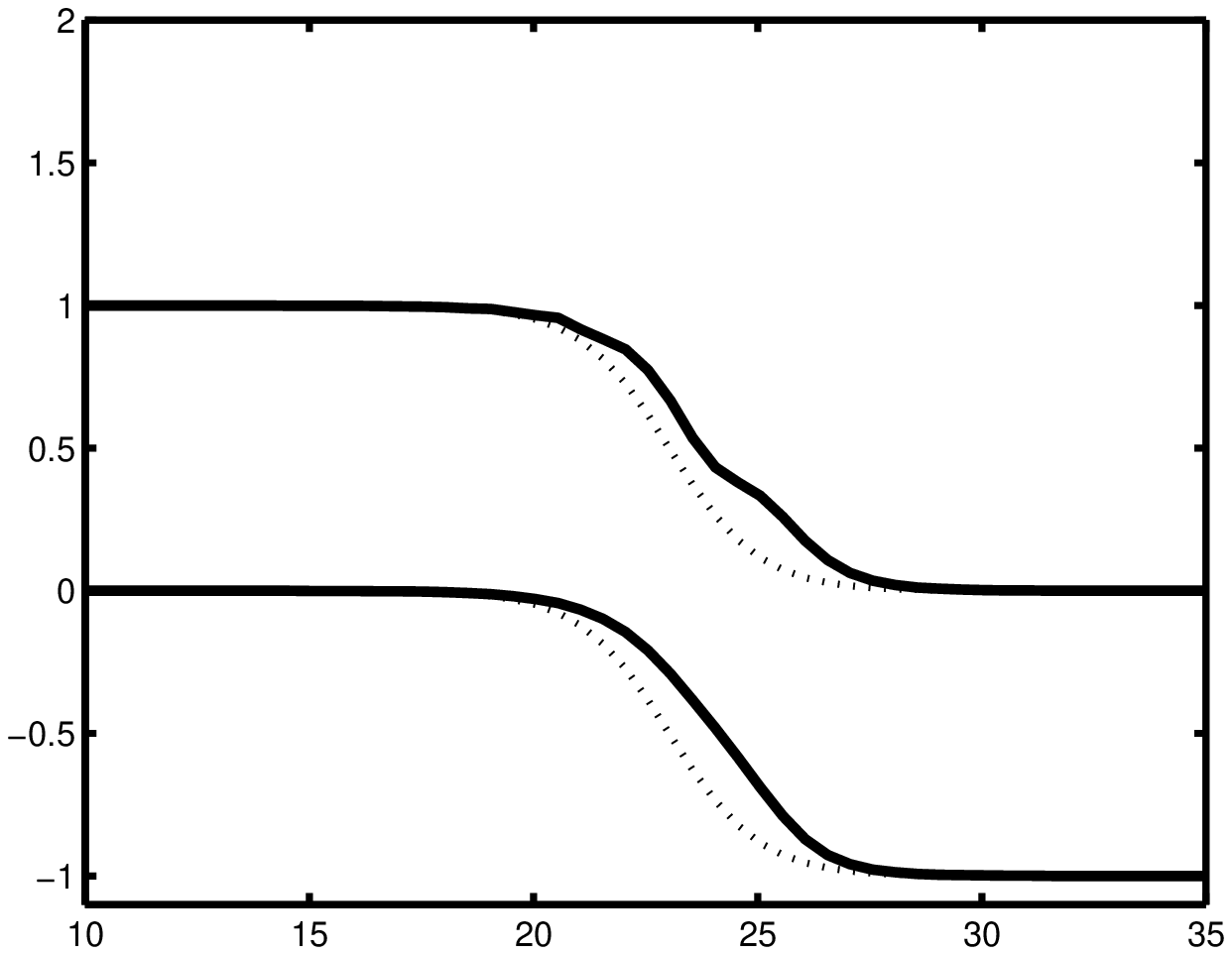} &
\includegraphics[width=5.75cm]{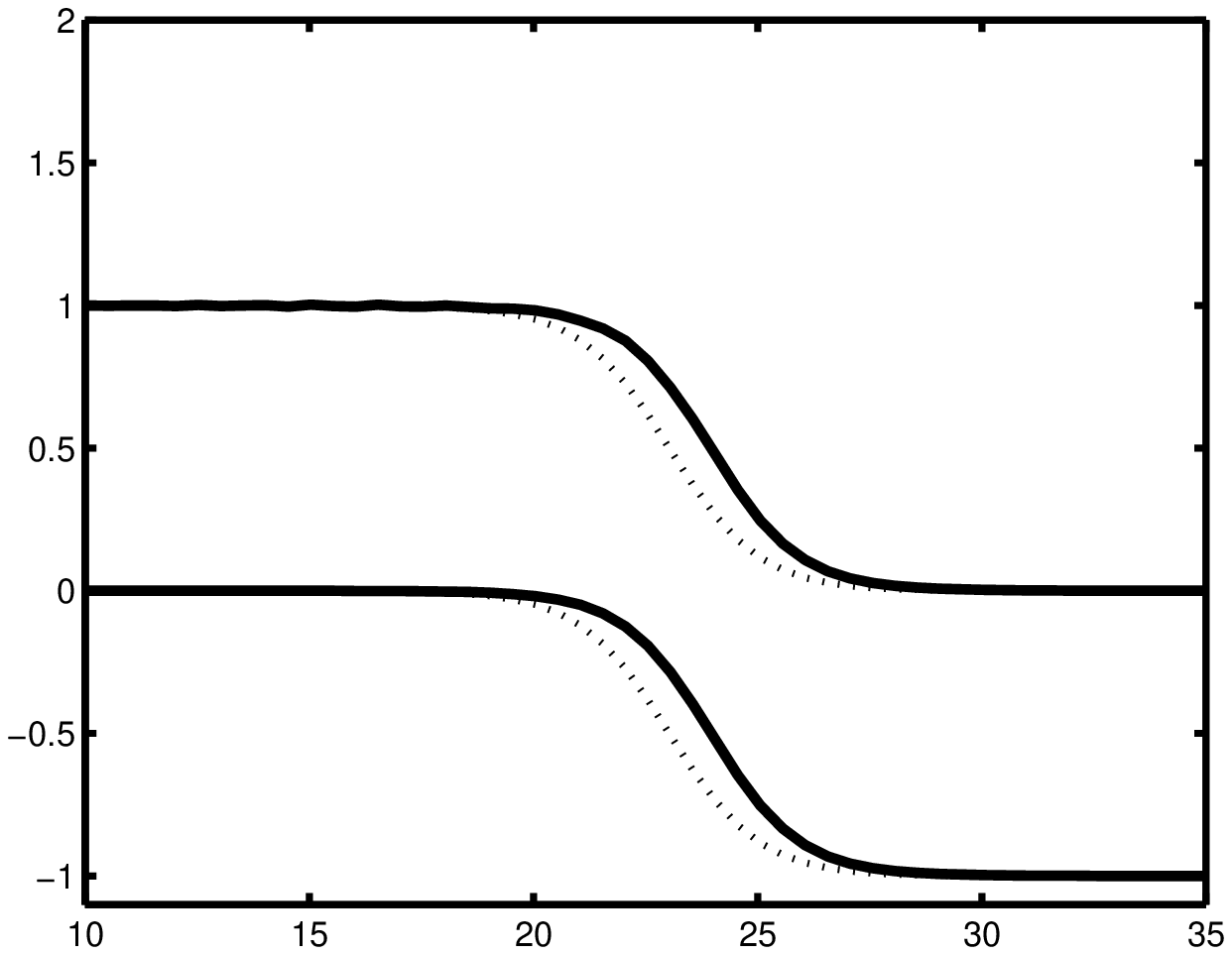}
\end{array}$
\end{center}
\caption{Snapshots of the evolution of a perturbed viscous shock wave 
solution generated by our extended Crank-Nicholson scheme (top curve:
$v$ against $x$; bottom curve: $u$ against $x$, with time $t$ fixed
and increasing from figure to figure). The parameters used were $\gamma = 1.4$ and $v_+ =9\times 10^{-6}$. This corresponds to a diatomic gas with Mach number $M \approx 2877$ (see Appendix \ref{mach}).  As expected, the wave converges to a translate of the original shock.}
\label{crank}
\end{figure}

 \[
\begin{split}
v_t+v_x-u_x&=0, \\
u_t+u_x-a\gamma v^{-\gamma -1}v_x &= \frac{u_{xx}}{v} - \frac{u_x v_x}{v^2}. 
\end{split}
\]
By implementing the Crank-Nicholson averaging, we obtain
\begin{align*}
&\frac{v_j^{n+1}-v_j^n}{\Delta t}+ \frac{1}{4 \Delta x}(v_{j+1}^{n+1} -v_{j-1}^{n+1} + v_{j+1}^n - v_{j-1}^n)\\
&\quad-\frac{1}{4 \Delta x}(u_{j+1}^{n+1} -u_{j-1}^{n+1} + u_{j+1}^n - u_{j-1}^n)=0
\end{align*}
and
\begin{align*}
&\frac{u_j^{n+1}-u_j^n}{\Delta t} + \frac{1}{4\Delta x}(u_{j+1}^{n+1} -u_{j-1}^{n+1} + u_{j+1}^n -u_{j-1}^n) \\
&\quad-\frac{a}{4\Delta x} \gamma (v_j^n)^{-\gamma-1}(v_{j+1}^{n+1}-v_{j-1}^{n+1} +v_{j+1}^n -v_{j-1}^n)\\
&\quad-\frac{1}{2(\Delta x)^2 v_j^n}(u_{j+1}^{n+1} -2u_j^{n+1} +u_{j-1}^{n+1} +u_{j+1}^n -2u_j^n + u_{j-1}^n)\\
&\quad+\frac{1}{16(\Delta x)^2 (v_j^n)^2} (u_{j+1}^{n+1} -u_{j-1}^{n+1}+ u_{j+1}^n -u_{j-1}^n)\\
&\quad \times (v_{j+1}^{n+1} -v_{j-1}^{n+1} + v_{j+1}^n -v_{j-1}^n)=0,
\end{align*}
where $n$ and $j$ are, respectively, the discretized temporal and spacial indices.

To cope with the nonlinearities, we use the Newton solver
\[
\begin{pmatrix}
U^{n+1}\\
V^{n+1}
\end{pmatrix}_{k+1} = 
\begin{pmatrix}
U^{n+1}\\
V^{n+1}
\end{pmatrix}_{k} - 
\begin{pmatrix}
D_uF & D_vF\\
D_uG & D_vG 
\end{pmatrix}^{-1}_{k},
\begin{pmatrix}
F(U^{n+1},V^{n+1})\\
G(U^{n+1},V^{n+1})
\end{pmatrix}_k
\]
where $F(U^{n+1},V^{n+1})$ and $G(U^{n+1},V^{n+1})$ are the above finite difference schemes, and $D_uF$, $D_vF$, $D_uG$, and $D_vG$ are their corresponding partial derivatives.  Hence, we use the previous time step as our initial guess in Newton's method and then iterate until convergence.

In Figure \ref{crank}, we see the evolution of a perturbed viscous shock.  We used a perturbation with a positive mass so that we could observe the convergence to a translate of the original profile, thus numerically demonstrating nonlinear stability.

\appendix
\section{Mach number for the $p$-system}\label{mach}

The Mach number is defined as
\[
M = \frac{u_+ - \sigma}{c_+},
\]
where $u_+$ is the downwind velocity, $\sigma$ is the shock speed, and $c_+$ is the downwind shock speed, all in Eulerian coordinates.  By considering the conservation of mass equation, we have $\rho_t + (\rho u)_x = 0$.  Hence, the jump condition is given by $\sigma [\rho] = [\rho u]$, which implies, in the original scaling for \eqref{psystem}, that
\[
\sigma = \frac{u_+ v_- - u_- v_+}{v_- - v_+}.
\]
Hence, 
\[
M = \frac{u_+ - \sigma}{c_+} = \frac{v_+(u_- - u_+)}{c_+(v_- - v_+)} = \frac{v_+ [u]}{c_+ [v]} = -s \frac{v_+}{c_+}
\]
or 
\[
M^2 = \left(\frac{u_+ - \sigma}{c_+}\right)^2 = \frac{s^2 v_+^2}{-p'(v_+)} = \frac{s^2 v_+^2}{\gamma a_0 v_+^{-\gamma-1}} = \frac{v_+^{\gamma + 3}}{\gamma} \frac{s^2}{a_0}.
\]
To express this in our scaling, which is given in \eqref{rescaled}, we need to swap the pluses and minuses.  Noting that $0<v_+<v_-=1$, we simplify to get
\[
M^2 = \frac{1}{\gamma v_+^\gamma} \frac{1-v_+^\gamma}{1-v_+} = \frac{1}{\gamma a}.
\]
Recalling that $a\sim v_+^\gamma$ as $v_+\to 0$, we
find that
\begin{equation}\label{vdep}
v_+\sim (\gamma M^2)^{-\frac 1\gamma}
\end{equation}
as $M\to \infty$.
In particular,
\begin{equation}\label{logvdep}
|\log v_+|\sim \gamma^{-1} (\log \gamma + 2\log M).
\end{equation}

\section{Profile bounds}\label{profbds}

Denote profile equation \eqref{profile} as
$ v' =H(v, v_+):= v(v-1 + a  (v^{-\gamma}-1)) $.

\begin{lemma}\label{xlem}
For $\gamma\ge 1$, $0\le x<1$, 
\begin{equation}
\label{xrel}
1\le \frac{1-x^\gamma}{1-x}\le \gamma.
\end{equation}
\end{lemma}
\begin{proof}
By convexity of $f(x)=x^\gamma$, the difference quotient \eqref{xrel}
is increasing in $x$, bounded above by $f'(1)=\gamma$, and below
by its value at $x=0$.
\end{proof}

\begin{proposition}
\label{lipcor}
For $\gamma\ge 1$ and $v_+\le v\le \frac{1}{6}$, $v_+\le \frac{1}{12}$,
\begin{equation}
\label{lipbd}
-\gamma(v-v_+)\le H(v, v_+)\le 
-\frac{3}{4} (v-v_+).
\end{equation}
For $\gamma\ge 1$, $v_+\le \frac{1}{4 \gamma}$,
and $\frac{3}{4}\le v\le v_-=1$,
\begin{equation}
\label{minuslipbd}
\frac{1}{2}(v-v_-)\le H(v, v_+)\le 
(v-v_-).
\end{equation}
\end{proposition}

\begin{proof}
By \eqref{RH}, 
\begin{align*}
H(v, v_+) &= v\Big((v-1) -
\frac{(v_+-1)(v^{-\gamma}-1)}{v_+^{-\gamma}-1}\Big)\\
&= v\Big((v-v_+) + \Big(\frac{1-v_+}{1-v_+^{\gamma}}\Big)
\Big(\big(\frac{v_+}{v}\big)^\gamma -1\Big)\Big)\\
&= (v-v_+)\Big(v - \Big(\frac{1-v_+}{1-v_+^{\gamma}}\Big)
\Big(\frac{1 - \big(\frac{v_+}{v}\big)^{\gamma}}{1 -
\big(\frac{v_+}{v}\big)}\Big)\Big).\\
\end{align*}
Applying \eqref{xrel} with 
$x=\frac{v_+}{v}$, we obtain \eqref{lipbd} from
\[
v- \gamma \le \frac{H(v,v_+)}{v-v_+}\le v-(1-v_+)
\]

Similarly, we may obtain \eqref{minuslipbd} by the calculation
$$
\frac{H(v,v_+)}{v-1}=
v\Big(
1-\Big(\frac{v_+}{v}\Big)^{\gamma}
\Big(\frac{1- v_+}{1-v_+^{\gamma}}\Big)
\Big(\frac{ 1-v^\gamma}{1-v}\Big)
\Big)
\ge v-\gamma v_+ \ge \frac{3}{4}- \frac{1}{4}.
$$
\end{proof}

\begin{corollary}
\label{expdecay}
For $\gamma\ge 1$, $0<v_+\le \frac{1}{12}$, and $\hat v(0):=v_+ + \frac{1}{12}$, the solution $\hat v$ of \eqref{profile} satisfies
\begin{subequations}
\label{decaybd}
\begin{align}
|\bV(x)-v_+|&\le \Big(\frac{1}{12}\Big)e^{-\frac{3x} {4}} \quad x\ge 0,\label{decaybd_1}\\
|\bV(x)-v_-|&\le 
\Big(\frac{1}{4}\Big)
e^{\frac{x+12}{2}} \quad x\le 0\label{decaybd_2}.
\end{align}
\end{subequations}
\end{corollary}

\begin{proof}
Observing that $(\hat v- v_+)(0)= \frac{1}{12}$, we obtain \eqref{decaybd_1} by Proposition \ref{lipcor} and the comparison principle for first-order scalar ODE.  Likewise, \eqref{decaybd_2} follows by \eqref{minuslipbd} together with the observation that, by convexity of $H$, $|H|$ is bounded below by estimates obtained at $\hat v=v_++\frac{1}{12}$ and $\hat v=\frac{3}{4}$ of $\frac{1}{16}$ and $\frac{1}{8}$, respectively, so that $\hat v$ traverses $[v_++\frac{1}{12}, \frac{3}{4}]$ over an $x$-interval of length $\le \frac{3/4}{1/16}= 12$.
\end{proof}

\begin{remark}
\label{regpert}
From Proposition \ref{lipcor} and Corollary \ref{expdecay}, we obtain the remarkable fact that, in the scaling we have chosen, $\bV$ decays up to first derivative to endstates $v_\pm$ at a uniform exponential rate {\it independent of shock strength}, despite the apparent singularity at $v=v_+$.
\end{remark}

\section{Initialization error} \label{initbds}

\begin{lemma}
For $A$ as in \eqref{evans_ode}, $|\cdot|$ the Euclidean ($\ell^2$) 
matrix operator norm,
\begin{subequations}
\label{Adecay}
\begin{align}
|A(x,\lambda)-A^+(\lambda)| &\le 
\Big( \frac{2|\lambda| + 1 + \gamma^2(\gamma-1) v_+^{-1}} {12}\Big)e^{- \frac{3x}{4}}, \quad x\ge 0,\label{Adecay_1}\\
|A(x,\lambda)-A^-(\lambda)| &\le 
\Big( \frac{2|\lambda| + 1 + 2\gamma^3(\gamma-1)}
{4}\Big)e^{\frac{x+12}{2}}, \quad x\le 0.
\label{Adecay_2}
\end{align}
\end{subequations}
\end{lemma}

\begin{proof}
By \eqref{evans_ode}, $|A(x,\lambda)-A_\pm(\lambda)|\le
2\lambda|\bV- v_\pm| + |f(\bV)-f(v_\pm)|$.
As computed in the proof of Lemma \ref{Hlem}, $f'(\bV)= -1 -a\gamma(\gamma-1)\bV^{-\gamma-1}.$
Applying \eqref{xrel} to the expression for $a$ in \eqref{RH},
we obtain $v_+^{\gamma}\le a \le \gamma v_+^{\gamma}$,
so that
\[
|f'(\bV)|\le 1 + \gamma^2(\gamma-1)\bV^{-1}
\le 1 + \gamma^2(\gamma-1)v_+^{-1},
\]
yielding \eqref{Adecay_1} by the Mean Value Theorem and \eqref{decaybd_1}.  Bound \eqref {Adecay_2} follows similarly.
\end{proof}

\begin{theorem}[Simplified Gap Lemma \cite{GZ,Br.1,Br.2}]\label{easygap}
Let $\widetilde{V}^+$ and $\widetilde{\mu}^+$ be a left eigenvector and
associated eigenvalue of $-A^+(\lambda)$ and suppose that
\begin{equation} \label{assume}
\begin{aligned}
|e^{(-A^+-\widetilde{\mu}^+)x}|&\le C_1 e^{-\hat\eta x} \quad x\le 0,\\
|(A-A^+)(x)|&\le C_2 e^{-\eta x} \quad x\ge 0,\\
\end{aligned}
\end{equation}
with $0\le \hat \eta <\eta$.
Then, there exists a solution $W=e^{\widetilde{\mu}^+ x}
\widetilde{V}(x, \lambda)$ of
\eqref{adjode} with
\begin{equation}\label{Vbd}
\frac{|\widetilde{V}(x,\lambda)- \widetilde{V}^+(\lambda)|}
{|\widetilde{V}^+(\lambda)|}
 \le \frac{C_1 C_2 e^{-\eta x}}{(\eta-\hat \eta) (1-\epsilon)} \quad x\ge L
\end{equation}
provided $(\eta-\hat \eta)^{-1} C_1 C_2 e^{-\eta L} \le \epsilon$.  
Similar estimates hold for solutions of \eqref{evans_ode} on $x\le 0$.
\end{theorem}

\begin{proof}
Writing $\widetilde{V}'=\widetilde{V}(-A^+ - \widetilde{\mu}^+) 
+ \widetilde{V}(-A + A^+)$ and imposing the limiting behavior 
$\widetilde{V}(+\infty, \lambda)=\widetilde{V}^+$, 
we obtain by Duhamel's Principle
\[
\widetilde{V}(x)=\widetilde{V}^+ -
\int_x^{+\infty} 
\widetilde{V}(y) 
e^{(-A^+ -\widetilde{\mu}^+)(x-y)}(A-A^+)
dy,
\]
from which the result follows by a straightforward Contraction Mapping argument using \eqref{assume}.  See \cite{GZ,Br.1,Br.2} for details.
\end{proof}

\begin{lemma}\label{C1}
For $\lambda$ satisfying \eqref{hfbounds}, $\gamma\in [1,3]$,
and $v_+\le \frac{1}{12}$,
\eqref{assume}(a) holds with $\eta= 1/2\gamma$, $\hat \eta=1/4\gamma$,
and $C_1=10^4$, and similarly for $x\le 0$.
\end{lemma}

\begin{proof}
Using the inverse Laplace transform representation
$$
e^{(-A^+-\widetilde{\mu}^+)x}=
\frac{1}{2\pi i}\oint_\Gamma e^{zt}(z+A^++\widetilde{\mu}^+)^{-1}dz,
$$
where $\Gamma$ is a contour enclosing the eigenvalues of
$(-A^+-\widetilde{\mu}^+)$ and distance $\hat \eta$ away,
and estimating the resolvent norm
$|(z+A^++\widetilde{\mu}^+)^{-1}|$ by Kramer's rule, we obtain
the stated crude bound, and similarly for $x\le 0$. 
\end{proof}

For error tolerance $\theta:=10^{-k}$ ($|\log \theta|\sim 2k$), define 
\begin{align*}
L_-(\theta) &:=2\Big( |\log 10^{-4}|+
|\log(2\gamma+7+ 2\gamma^3(\gamma-1))|+  |\log \theta| \Big)
+ 12,\\
L_+(\theta, v_+)&:=
\frac{4}{3}\Big( |\log 10^{-4}|+
|\log( 2\gamma+7+ \gamma^2(\gamma-1) v_+^{-1})| + |\log \theta|
\Big).
\end{align*}

\begin{corollary}\label{initcor}
For $\lambda$ satisfying \eqref{hfbounds}, we have relative error bounds
\begin{equation}
\label{relerr}
\frac{|W_1^-(-L_-,\lambda)- V^-e^{\mu^- x}|}{|V^-e^{\mu^- x}|},
\quad
\frac{|\widetilde{W}_1^+(L_+,\lambda)- \widetilde{V}^+e^{\widetilde{\mu}^+ x}|} {|\widetilde{V}^+ e^{\widetilde{\mu}^+ x}|} \le \theta.
\end{equation}
\end{corollary}

\begin{proof}
This follows by 
\eqref{Vbd} with \eqref{Adecay}, Lemma \ref{C1}, 
and \eqref{hfbounds}.
\end{proof}

\begin{remark}\label{Lest}
Combining \eqref{relerr} with \eqref{logvdep},
we find that
\[
L_+\sim \frac{4}{3}(  2 \log M + 4+ |\log 10^{-4}|+ |\log \theta|)
\]
suffices to obtain relative initialization error less than $\theta$
for $\gamma\in [1,3]$ and $\lambda$ in the 
computational region \eqref{hfbounds}.
For $\theta=10^{-3}$, we thus obtain $\theta$ tolerance up to $M=3,000 \sim 10^3$ for $L_+\sim 40$. 
This is conservative, as we have made little effort to optimize bounds, but still within the realm of our experiments.  Our numerical convergence studies indicate that $L_\pm=18$
in fact suffices for $10^{-3}$ accuracy.
\end{remark}

\def\ocirc#1{\ifmmode\setbox0=\hbox{$#1$}\dimen0=\ht0 \advance\dimen0
  by1pt\rlap{\hbox to\wd0{\hss\raise\dimen0
  \hbox{\hskip.2em$\scriptscriptstyle\circ$}\hss}}#1\else {\accent"17 #1}\fi}
  \def\cprime{$'$}


\begin{thebibliography}{10}

\bibitem{AGJ}
J.~Alexander, R.~Gardner, and C.~Jones.
\newblock A topological invariant arising in the stability analysis of
  travelling waves.
\newblock {\em J. Reine Angew. Math.}, 410:167--212, 1990.

\bibitem{AS}
J.~C. Alexander and R.~Sachs.
\newblock Linear instability of solitary waves of a {B}oussinesq-type equation:
  a computer assisted computation.
\newblock {\em Nonlinear World}, 2(4):471--507, 1995.

\bibitem{Ba}
G.~K. Batchelor.
\newblock {\em An introduction to fluid dynamics}.
\newblock Cambridge Mathematical Library. Cambridge University Press,
  Cambridge, paperback edition, 1999.

\bibitem{BSZ}
S.~Benzoni-Gavage, D.~Serre, and K.~Zumbrun.
\newblock Alternate {E}vans functions and viscous shock waves.
\newblock {\em SIAM J. Math. Anal.}, 32(5):929--962 (electronic), 2001.

\bibitem{BDG}
T.~J. Bridges, G.~Derks, and G.~Gottwald.
\newblock Stability and instability of solitary waves of the fifth-order
  {K}d{V} equation: a numerical framework.
\newblock {\em Phys. D}, 172(1-4):190--216, 2002.

\bibitem{Br.1}
L.~Q. Brin.
\newblock {\em Numerical testing of the stability of viscous shock waves.}
\newblock PhD thesis, Indiana University, Bloomington, 1998.

\bibitem{Br.2}
L.~Q. Brin.
\newblock Numerical testing of the stability of viscous shock waves.
\newblock {\em Math. Comp.}, 70(235):1071--1088, 2001.

\bibitem{BrZ}
L.~Q. Brin and K.~Zumbrun.
\newblock Analytically varying eigenvectors and the stability of viscous shock
  waves.
\newblock {\em Mat. Contemp.}, 22:19--32, 2002.
\newblock Seventh Workshop on Partial Differential Equations, Part I (Rio de
  Janeiro, 2001).

\bibitem{E.1}
J.~W. Evans.
\newblock Nerve axon equations. {I}. {L}inear approximations.
\newblock {\em Indiana Univ. Math. J.}, 21:877--885, 1971/72.

\bibitem{E.2}
J.~W. Evans.
\newblock Nerve axon equations. {II}. {S}tability at rest.
\newblock {\em Indiana Univ. Math. J.}, 22:75--90, 1972/73.

\bibitem{E.3}
J.~W. Evans.
\newblock Nerve axon equations. {III}. {S}tability of the nerve impulse.
\newblock {\em Indiana Univ. Math. J.}, 22:577--593, 1972/73.

\bibitem{E.4}
J.~W. Evans.
\newblock Nerve axon equations. {IV}. {T}he stable and the unstable impulse.
\newblock {\em Indiana Univ. Math. J.}, 24(12):1169--1190, 1974/75.

\bibitem{EF}
J.~W. Evans and J.~A. Feroe.
\newblock Traveling waves of infinitely many pulses in nerve equations.
\newblock {\em Math. Biosci.}, 37:23--50, 1977.

\bibitem{GJ}
R.~Gardner and C.~K. R.~T. Jones.
\newblock A stability index for steady state solutions of boundary value
  problems for parabolic systems.
\newblock {\em J. Differential Equations}, 91(2):181--203, 1991.

\bibitem{GZ}
R.~A. Gardner and K.~Zumbrun.
\newblock The gap lemma and geometric criteria for instability of viscous shock
  profiles.
\newblock {\em Comm. Pure Appl. Math.}, 51(7):797--855, 1998.

\bibitem{Go.1}
J.~Goodman.
\newblock Nonlinear asymptotic stability of viscous shock profiles for
  conservation laws.
\newblock {\em Arch. Rational Mech. Anal.}, 95(4):325--344, 1986.

\bibitem{HZ}
P.~Howard and K.~Zumbrun.
\newblock Pointwise estimates and stability for dispersive-diffusive shock
  waves.
\newblock {\em Arch. Ration. Mech. Anal.}, 155(2):85--169, 2000.

\bibitem{HSZ}
J.~Humpherys, B.~Sandstede, and K.~Zumbrun.
\newblock Efficient computation of analytic bases in {E}vans function analysis
  of large systems.
\newblock {\em Numer. Math.}, 103(4):631--642, 2006.

\bibitem{HuZ.2}
J.~Humpherys and K.~Zumbrun.
\newblock An efficient shooting algorithm for evans function calculations in
  large systems.
\newblock {\em Physica D}, 220(2):116--126, 2006.

\bibitem{Kato}
T.~Kato.
\newblock {\em Perturbation theory for linear operators}.
\newblock Classics in Mathematics. Springer-Verlag, Berlin, 1995.
\newblock Reprint of the 1980 edition.

\bibitem{MZ.1}
C.~Mascia and K.~Zumbrun.
\newblock Pointwise {G}reen function bounds for shock profiles of systems with
  real viscosity.
\newblock {\em Arch. Ration. Mech. Anal.}, 169(3):177--263, 2003.

\bibitem{MZ.3}
C.~Mascia and K.~Zumbrun.
\newblock Stability of large-amplitude viscous shock profiles of
  hyperbolic-parabolic systems.
\newblock {\em Arch. Ration. Mech. Anal.}, 172(1):93--131, 2004.

\bibitem{MZ.2}
C.~Mascia and K.~Zumbrun.
\newblock Stability of small-amplitude shock profiles of symmetric
  hyperbolic-parabolic systems.
\newblock {\em Comm. Pure Appl. Math.}, 57(7):841--876, 2004.

\bibitem{MN}
A.~Matsumura and K.~Nishihara.
\newblock On the stability of travelling wave solutions of a one-dimensional
  model system for compressible viscous gas.
\newblock {\em Japan J. Appl. Math.}, 2(1):17--25, 1985.

\bibitem{PSW}
R.~L. Pego, P.~Smereka, and M.~I. Weinstein.
\newblock Oscillatory instability of traveling waves for a {K}d{V}-{B}urgers
  equation.
\newblock {\em Phys. D}, 67(1-3):45--65, 1993.

\bibitem{PW}
R.~L. Pego and M.~I. Weinstein.
\newblock Eigenvalues, and instabilities of solitary waves.
\newblock {\em Philos. Trans. Roy. Soc. London Ser. A}, 340(1656):47--94, 1992.

\bibitem{Sm}
J.~Smoller.
\newblock {\em Shock waves and reaction-diffusion equations}.
\newblock Springer-Verlag, New York, second edition, 1994.

\bibitem{TZ.1}
B.~Texier and K.~Zumbrun.
\newblock Relative {P}oincar\'e-{H}opf bifurcation and galloping instability of
  traveling waves.
\newblock {\em Methods Appl. Anal.}, 12(4):349--380, 2005.

\bibitem{TZ.2}
B.~Texier and K.~Zumbrun.
\newblock Galloping instability of viscous shock waves.
\newblock Preprint, 2006.

\bibitem{TZ.3}
B.~Texier and K.~Zumbrun.
\newblock Hopf bifurcation of viscous shock waves in compressible gas- and
  magnetohydrodynamics.
\newblock Preprint, 2006.

\bibitem{Z.2}
K.~Zumbrun.
\newblock Multidimensional stability of planar viscous shock waves.
\newblock Lecture notes, TMR summer school, 1999.

\bibitem{ZH}
K.~Zumbrun and P.~Howard.
\newblock Pointwise semigroup methods and stability of viscous shock waves.
\newblock {\em Indiana Univ. Math. J.}, 47(3):741--871, 1998.

\end{thebibliography}
\end{document}